\documentclass{svproc}
\usepackage{url}

\usepackage{subfigure}
\usepackage{url}
\usepackage{amsmath,amsfonts,amssymb}
\usepackage{array,multirow}
\usepackage{graphics}
\usepackage{graphicx}
\usepackage{float}
\usepackage{caption}
\usepackage{verbatim}
\usepackage{color}
\usepackage{rotating}
\usepackage{tabularx}
\usepackage{lscape}
\usepackage{url}
\usepackage{algorithmic}
\usepackage{booktabs}
\usepackage{rotating}
\usepackage{sidecap}
\usepackage[boxed]{algorithm2e}
\usepackage{hyperref}

\begin{document}
\mainmatter              % start of a contribution
\title{Optimizing the nozzle path in the 3D Printing Process}
\titlerunning{Optimizing the nozzle path in the 3D Printing Process}  % abbreviated title (for running head)
\author{Manuel Iori\inst{1} \and Stefano Novellani\inst{1}}
\authorrunning{Iori and Novellani} % abbreviated author list (for running head)
%\tocauthor{Manuel Iori and Stefano Novellani}
%
\institute{Univesity of Modena and Reggio Emilia,\\ DISMI, Via Amendola 2, 42123 Reggio Emilia, Italy,\\
\email{manule.iori;stefano.novellani@unmimore.it}}

\maketitle              % typeset the title of the contribution

\begin{abstract}
In this paper, we define the 3D printing routing problem, the problem of finding the optimal path of the nozzle in a fused deposition modeling 3D printing system, so as to minimize the time required to create on object. We formally model the problem with an integer linear programming formulation and then solve it via heuristic algorithms. We test the algorithms on a set of large-size real-life instances, comparing them with one of the most widely used open source software for the problem. We show that large time reductions can be obtained. We finally propose a set of interesting directions for future research.   \\
\end{abstract}

\section{Introduction}\label{sec:intro3d}
{\em Three-dimensional Printing} (3DP), also known as {\em Additive Manufacturing}, is a production technology that allows to construct three-dimensional objects with complex shapes by adding one or more materials layer by layer.
Since its rise in the 80s, different techniques for 3DP have been developed, including {\em Stereolithography}, that makes use of liquid polymers hardened by UV beams layer by layer, {\em Selective Laser Sintering}, where a laser melts the powder, layer by layer, in a selective way, and {\em Fused Deposition Modeling} (FDM), also known as {\em Fused Filament Fabrication}.
FDM involves extruding thermoplastic material (such as ABS plastic or polycarbonate) through heated nozzles: the filament is pulled by a feed roller, it is heated to a temperature above the material's melting point, and it is then extruded by the nozzle. Once a layer is completed, the platform where the material has been deposited is lowered (or the nozzle is raised) to the next layer level, and a new layer is added above the previous one.
%A supporting material can also be used in the process when it is needed for structural issues, especially in case of overhanging geometries
From now on, when mentioning 3DP we will refer to the FDM technique.

In this work, we focus on the optimization of the 3DP process, in particular we concentrate on the optimization of the path that the nozzle has to follow to produce the three-dimensional object. {The objective is to minimize the printing time while respecting {given} printing constraints.
{The large amount of time needed to produce objects with 3DP is one of the main disadvantages of this technology, representing a limit to its use for mass production,
so minimizing its printing time is a very relevant objective.}
Before introducing models and heuristic algorithms, we describe the 3DP process. We focus on the aspects that are related with the optimization pursued in this paper, and refer the reader interested in structural properties of the materials and technical details of the production process to, e.g., Campbell et al. \cite{campbell2011could} and Akella \cite{akella2012professor}.

%\subsection{The 3DP Process}
The manufacturing process of a product through 3DP follows the typical scheme depicted in Figure \ref{fig:3Dprocess}. The process starts with the definition of the 3D model, which can be created by a computer-aided design (CAD) software or by directly scanning an object that we want to reproduce. The second step is to reproduce the three-dimensional object in a polygonal mesh, normally made of triangles. The typical file used to report a polygonal mesh is in a .stl format. When the model is represented as a mesh, the corresponding .stl file is passed to a slicing software that divides the mesh in layers. The following step is the most interesting to our work and it is to decide the path to be followed by the tool for reproducing the three-dimensional object. This is usually obtained by dedicated software/freeware by the use of quick heuristic algorithms (see, e.g., Cura \cite{cura} and Slic3r \cite{slicer}). The solution obtained is expressed {as} a sequence of commands in a g-code format, {which is then} passed to the printer in order to produce the final 3D object (see, e.g., Campbell et al. \cite{campbell2011could}).
\begin{figure}[h]
\centering
\includegraphics[width=1\textwidth]{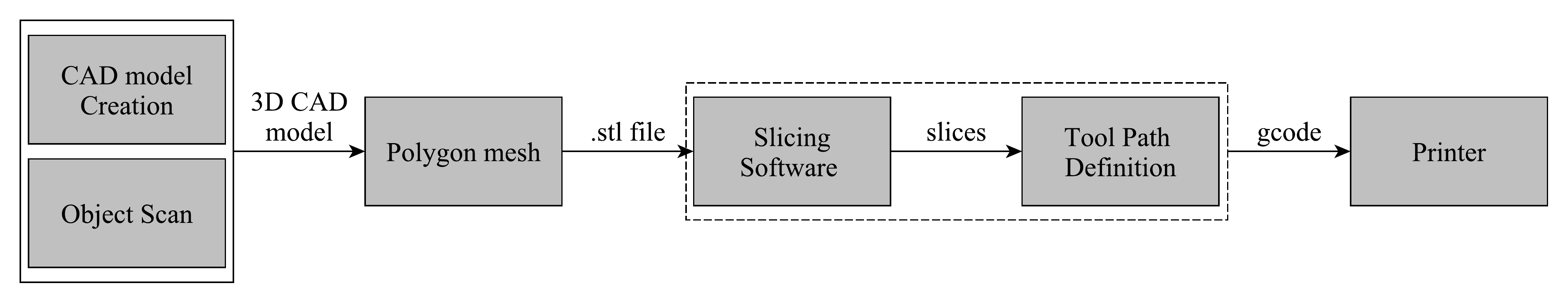}
\caption{Scheme of the 3DP process. %(ad un certo punto fare un disegnino di tutto ci\'o)
}
\label{fig:3Dprocess}
\end{figure}

In the 3DP process, there are {two relevant} aspects that can be optimized. The first one is about the quality of the printed object and {in particular its surface}.
The second one {refers to the printing time} that can be minimized, and it is strictly connected to the tool path. In this paper, we account for this second aspect by focusing in minimizing the distance traveled by the nozzle during the 3DP process, {thus indirectly} minimizing the time of the print.

{In the following}, we present a formal description of our 3DP optimization problem, which can be seen as a generalization of the well-known rural postman problem (see, e.g., Eiselt et al. \cite{eiselt1995arc}). We then {propose a brief} literature review and a mathematical model (Sections \ref{sec:liter} and \ref{sec:form}).} To solve the problem, we introduce some heuristic algorithms and produce a collection of instances derived from three-dimensional objects on which we perform computational tests (Sections \ref{sec:heur} and \ref{sec:comp}). Last but not least, we list a set of features and future research direction to be followed in the field that represents the intersection between optimization and 3DP.

%The paper is organized as follows. {In Section \ref{sec:liter} we will define the class of the Rural Postman Problems and produce an exhaustive literature review. In Section \ref{sec:probdesc} we define the problem we are presenting, while in Section \ref{sec:form} we give a formulation for the problem. In Section \ref{sec:heur} we present some heuristic algorithms to solve the problem and in Section \ref{sec:comp} we show the computational results of the presented algorithms on a set of newly collected instances.}

\section{The Rural Postman Problem} \label{sec:liter}
The problem of defining the tool path in a 3DP process is a generalization of the {\em Rural Postman Problem} (RPP), that is a particular version of the more general class of {\em Arc Routing Problems} (see, e.g., Corber\'an and Laporte \cite{corberan2014arc}).
Given a multigraph $G = (V,A)$, where $V$ is the set of vertices and $A$ the set of arcs (that can be directed, indirected, or both), {each associated with a non-negative cost, and given a subset $R\subseteq A$ of required arcs, the RPP is the problem of finding a minimum-cost closed path that traverses all arcs of $R$ at least once}. The RPP is NP-hard if the set of required arcs is not strongly connected (Laporte \cite{laporte1997modeling}).
%The RPP represents a generalization of the {\em Chinese Postman Problem} (CPP), that is the problem of finding the minimum cost path that visits at least once all the arcs. The CPP can be solved in polynomial time, while the RPP is an NP-hard problem if the set of required arcs is not strongly connected (Laporte \cite{laporte1997modeling}).

The RRP arises in many real-world situations such as  snow plowing,  garbage collection, etc. The reader can find several contributions related to these applications in Eiselt et al. \cite{eiselt1995arc}. More recent applications are referred to the control of plotting and drilling machines (see, e.g., Groetschel \cite{grotschel1991optimal}) and to the optimization of laser-plotter beam movements (see, e.g., Ghiani and Improta \cite{ghiani2001laser}).

Many polyhedral studies and exact algorithms, mostly cutting planes and branch-and-cut algorithms, have been proposed in the last years for the RRP and its generalizations.
Laporte and Ghiani \cite{ghiani2000branch} solved the {\em undirected RPP} by means of a branch-and-cut algorithm based on a formulation making use of only binary variables. Fern\'andez et al. \cite{fernandez2003undirected} proposed a {novel formulation} and solved it with an algorithm that derives the lower bounds via cutting planes and the upper bounds from heuristic algorithms.
According to the survey by Corber\'an and Prins \cite{corberan2010recent}, the largest RPP instances that have been solved to optimality are characterized by up to 1000 vertices, 3080 edges and 204 required components and are solved in about one hour {of computing time on a standard PC by using the branch-and-cut described in Corber\'an et al. \cite{corberan2007branch}} for the {\em windy general routing problem} (which contains the RPP).
%In Corber\'an and Sanchis \cite{corberan1994polyhedral,corberan1998general}, the authors proposed a polyhedral approach and defined some facets for the RPP.
%Reinelt and Theis  \cite{reinelt2006note} improved the results of Laporte and Ghiani.

Probably, the most famous heuristic algorithm for the RPP is the one by Frederikson \cite{frederickson1979approximation}, which is built on the well-known heuristic for the {\em Traveling Salesman Problem} (TSP) by Christofides \cite{christofides1976worst}, and, similarly, presents a worst case bound of 3/2. Hertz et al. \cite{hertz1999improvement} proposed some improving techniques to be applied to the Frederickson algorithm and to a newly developed heuristic. Corber\'an et al. \cite{corberan2000heuristics} presented two heuristic approaches to solve the {\em Mixed RPP} (MRPP), {which is the RPP defined on a graph made of both edges and arcs}. Groves and Vuuren \cite{groves2005efficient} presented an effective local search framework for the undirected RPP based on local searches for the TSP. They solved heuristically very large instances with up to about 5000 vertices and 30000 edges. Ghiani et al. \cite{ghiani2006constructive} proposed a constructive heuristic for the {\em undirected RPP} with a local post-optimization. The procedure is competitive with the Frederickson one. Holmberg \cite{holmberg2010heuristics} proposed some heuristics for the RPP built on the Frederikson heuristic by changing the order of the algorithmic components. The author includes also two post-processing heuristics to improve the solution.

We report now a set of variants of the RPP that are interesting for 3DP. {Dror and Langevin \cite{dror1997generalized} introduced the {\em directed clustered RPP}, a generalization of the RPP in which each subset of arcs to be served has to be completed before serving another subset. These subsets do not have to be visited in a given order. An application of this problem is the postal delivery, where the delivery is divided into subsets of clients. The authors proposed an enumerative solution approach by transforming the directed clustered RPP into a {\em Generalized TSP}, that is a modification of the TSP in which nodes are partitioned into clusters and exactly one node from each cluster is visited in a cycle. The largest instance they solved is made of 581 vertices, 770 arcs, and 44 required arcs.}

Letchford and Eglese \cite{letchford1998rural} studied the RPP with deadline classes, a RPP with time windows {where a maximum time has to be respected for each vertex}. They proposed a formulation and a cutting plane algorithm. Corber\'an et al. \cite{corberan2002rural} consider the mixed RPP with turn penalties. They associated a penalty to every turn and took into account the existence of forbidden turns. They transformed the problem into an {\em Asymmetric TSP} (ATSP) and solved it by using exact and heuristic algorithms. Ghiani and Improta \cite{ghiani2001laser} considered the problem of drawing and decorating metal surfaces by means of a laser plotter that works as a 2D printer, minimizing the total length of the non-drawing moves. The so-called {\em Laser-plotter beam Routing Problem} (LRP) is modeled as an RPP with additional constraints {that aim at minimizing the spot created on the surface every time a new line is drawn and at reducing the number of the shutter openings to reduce maintenance costs. They solved the LRP by transforming it into an RPP.} By using the branch-and-cut presented in Ghiani and Laporte \cite{ghiani2000branch}, instances with up to 225 vertices were solved. Moreira et al. \cite{moreira2007heuristics} described the problem of determining the shortest path for the cutting of given pieces when manufacturing with high-precision tools. The problem is represented as a particular RPP where non-cutting movements are allowed only after the complete cut of a piece. The problem is called dynamic RPP and solved by means of a heuristic algorithm. Orazi et al. \cite{orazi2015} optimized the path of the scanning head in laser texturing of free form surfaces. They modeled the problem as a TSP and solved it by using the heuristic methods provided by the Concorde software.

The following recent works have considered optimization in the 3D process. Panesar et al. \cite{panesar2014design} designed an optimization strategy for the definition of the internal configuration for structural purpose.  Zhang et al. \cite{zhang2015} and Wang et al. \cite{wang2016} studied the best printing direction to diminish the possible flaws and improve the surface quality given by the external support needed when printing an object.

\section{The 3D Printing Routing Problem} \label{sec:form}

We call {\em 3D printing Routing Problem} (3DRP) the generalization of the RPP that is related to the tool path definition in a 3DP process.
The 3DRP represents the {\it basic problem} of optimizing the tool path of a 3DP process, but some additional characteristics {might arise across the 3DP process. We provide} a brief description of those features in Section \ref{future}.

{In the 3DP tool path, the nozzle must start from its starting position (vertex 0) and return to the same position when the print is terminated.
The graph is divided into clusters to be visited in a sequence, {each representing a layer to be printed}. On each layer, the set of {\it compulsory edges} represents the set of edges on which the nozzle is required to depose material, while the other edges are called {\it optional edges}.
No layer can be started before terminating the previous one, and the nozzle cannot move back to previous layers. Then, being that only optional edges are present between two layers, one and only one optional edge can be used to move from one layer to the next one. Moreover, no edge between two non-subsequents layers exists.}

{Formally}, the problem is defined on an undirected graph $G = (V, E)$, where $V$ is the set of vertices and $E$ is the set of edges. Moreover, it is defined a subset $R \subseteq E$ of edges that must be visited at least once. In our case, {$R$ contains} the compulsory edges, where the material is extruded (printed) the first time they are used, {whereas $E \setminus R$ contains} the set of optional edges.
Let us define {\em layer} as a set of vertices characterized by the same altitude.
We call $\mathcal{L}$ the set of all layers, where a layer $L_h \in \mathcal{L}$ is one layer of print, so $\mathcal L = \cup_{h = 1}^{|\mathcal L|} L_h$, and $L_h \cap L_k = \emptyset, h \ne k, L_h, L_k \in \mathcal{L}$. The graph $G$ is clustered in several subgraphs, {each induced by a layer plus the starting and ending points} of print, vertex 0, and all the edges connecting these components. Thus, $V = \mathcal{L} \cup \{0\}$.
%
%\begin{align}
%\min  \sum_{(i,j)\in E} c_{ij} x_{ij} \label{eq:fobb}\\
%x_{ij} \ge 1 & & (i,j) \in R \label{eq:compulsory}\\
%\sum_{j \in V : j > i} x_{ij} + \sum_{j \in V : j < i} x_{ji} = 2z && i \in V  \label{eq:flow}\\
%\sum_{i \in S} \sum_{j \in \bar S} x_{ij} \ge 1 & & S \subseteq L_h, \bar S = L_{h+1} \cup L_{h} \setminus S, L_h \subseteq \mathcal{L} \setminus L_{|\mathcal L|} \label{eq:connectivity}\\
%\sum_{i \in S} \sum_{j \in \bar S} x_{ij} \ge 1 & & S \subseteq L_{|\mathcal L|}, \bar S = \{0\} \cup L_{|\mathcal L|}\setminus S \label{eq:connectivity2}\\
%\sum_{j \in L_1} x_{0j} = 1 \label{eq:outdepot}\\
%\sum_{j \in L_{|\mathcal L|}} x_{j0} = 1 \label{eq:indepot}\\
%\sum_{(i,j):i \in L_h, j \in L_{h+1}} x_{ij} = 1 & & L_h \subseteq \mathcal{L} \setminus L_{|\mathcal L|} \label{eq:infralayer}\\
%x_{ij} \ge 0 & & (i,j) \in E  \label{eq:xpos}\\
%z \ge 0 &&  \label{eq:zpos}
%\end{align}
%
%

Let variable $x_{ij}$ give the number of times edge  $(i,j) \in E \setminus R$ is traveled, and the number of times edge $(i,j) \in R$ is traveled in addition to the first time. Let $c_{ij}$ be the shortest distance between $i$ and $j$, $(i,j) \in E$. The 3DRP aims at minimizing the path of the nozzle starting from vertex 0, traveling all the compulsory edges at least once, respecting the order of the layers, and returning to the vertex 0.
{Let us define
$\Omega = \{S \subseteq L_h, \bar S = L_{h+1} \cup L_{h} \setminus S, \text{for } h=1,\dots,|\mathcal L|-1 \} \cup \{S \subseteq L_{|\mathcal L|}, \bar S = \{0\} \cup L_{|\mathcal L|}\setminus S\}$. Then, for each $S \in \Omega$, let us also define
$R^1(S) = \{(i,j) \in R: i \in S, j \in \bar S, j > i\}$,
$R^2(S) = \{(i,j) \in R: i \in S, j \in \bar S, j < i\}$,
$E^1(S) = \{(i,j) \in E \setminus R: i \in S, j \in \bar S, j > i\}$,
$E^2(S) = \{(i,j) \in E \setminus R: i \in S, j \in \bar S, j < i\}$. The basic 3DRP can be modeled as}
{\small
\begin{align}
\allowdisplaybreaks
&\min  \sum_{(i,j)\in R} c_{ij}(1 + x_{ij}) +  \hspace{-0.3cm} \sum_{(i,j)\in E \setminus R} c_{ij} x_{ij} \label{eq:fobb}\\
& \sum_{j \in R : j > i} \hspace{-0.1cm} (1 + x_{ij}) + \hspace{-0.3cm} \sum_{j \in R : j < i} \hspace{-0.1cm} (1 + x_{ji}) + \hspace{-0.3cm} \sum_{j \in V : j > i} \hspace{-0.1cm} x_{ij} + \hspace{-0.3cm} \sum_{j \in V : j < i} \hspace{-0.1cm} x_{ji} = 2z_i, & \, i \in V  \label{eq:flow}\\
%&\sum_{\substack{(i,j) \in R : i \in S, \\ j \in \bar S, j > i}} \hspace{-0.3cm} (1 + x_{ij}) + \hspace{-0.3cm} \sum_{\substack{(i,j) \in R : i \in S, \\ j \in \bar S, j < i}} \hspace{-0.3cm}  (1 + x_{ij}) + \hspace{-0.3cm} \sum_{\substack{(i,j) \in E \setminus  R : i \in S, \\ j \in \bar S, j > i}} \hspace{-0.5cm}  x_{ij} + \hspace{-0.3cm} \sum_{\substack{(i,j) \in E \setminus R : i \in S,\\ j \in \bar S, j > i}} \hspace{-0.5cm}x_{ij}\ge 1  \label{eq:connectivity}\\
&\sum_{{(i,j) \in R^1(S)}} \hspace{-0.3cm} (1 + x_{ij}) + \hspace{-0.5cm} \sum_{{(i,j) \in R^2(S)}} \hspace{-0.3cm}  (1 + x_{ij}) + \hspace{-0.5cm} \sum_{{(i,j) \in E^1(S)}} \hspace{-0.3cm}  x_{ij} + \hspace{-0.5cm} \sum_{{(i,j) \in E^2(S)}} \hspace{-0.3cm}x_{ij} \ge 1 & \, S \subseteq \Omega \label{eq:connectivity}\\
%& (S \subseteq L_h, \bar S = L_{h+1} \cup L_{h} \setminus S,  h=1,\dots,|\mathcal L|-1)  \lor
%(S \subseteq L_{|\mathcal L|}, \bar S = \{0\} \cup L_{|\mathcal L|}\setminus S) \label{eq:connectivity}  \\
&\sum_{j \in L_1} x_{0j} = \sum_{j \in L_{|\mathcal L|}} x_{j0} = 1 \label{eq:outdepot}\\
%\sum_{j \in L_{|\mathcal L|}} x_{j0} = 1 \label{eq:indepot}\\
&\sum_{(i,j):i \in L_h, j \in L_{h+1}}  \hspace{-0.3cm} x_{ij} = 1 & \,   h=1,...,|\mathcal L|-1 \label{eq:infralayer}\\
&x_{ij} \in \mathbb{N}_0 & \, (i,j) \in E  \label{eq:xpos}\\
&z_i \in \mathbb{N}_0 & \, i \in V \label{eq:zpos}
\end{align}
}
The objective function \eqref{eq:fobb} aims at minimizing the total costs, i.e., distances. Constraints \eqref{eq:flow} state that the number of edges incident to a vertex are even or zero. In \eqref{eq:connectivity}, we impose the connectivity of the solution. Constraints \eqref{eq:outdepot} state that one and only one edge must leave and return to the vertex 0. Constraint \eqref{eq:infralayer} imposes that one and only one edge can be used between one layer and the next one. The variables are defined as integer and non-negative in  \eqref{eq:xpos} and \eqref{eq:zpos}.

\section{Heuristic Algorithms} \label{sec:heur}
The 3DP process {typically} involves from thousands to millions of edges to be printed and traveled by the nozzle. Thus, real 3DRP instances are impossible to solve exactly, because of the huge number of variables and constraints needed. For this reason, we developed four heuristic algorithms.
%To do so, we organized the input in couples of vertices that defines the edges to be traveled and printed on each layer (the set $R$), named as compulsory edges. The optional edges ($E \setminus R$) and the compulsory edges used more than once are the {\it non-printing edges}, and are computed and inserted just when needed during the construction of the solution.
%The input is obtained by analyzing the g-codes created by Cura \cite{cura}, one of the most advanced 3DP software.

{\bf Closest 3DP.}
The first algorithm we present to solve the 3DRP starts from vertex 0 and moves from one edge to the closest unused compulsory edge (we enter in the edge by the closest vertex with respect to the current one and then we travel the edge to arrive in its second vertex, which will be the new current vertex). If the second vertex of the current edge is also the first vertex of the next edge, we link them directly, otherwise we insert a non-printing edge. When the compulsory edges of the current layer have been traversed at least once, we move to the next layer by using the cheapest edge. We repeat this until the last layer is completed, and then connect the last visited vertex to vertex 0 to end the algorithm.

\begin{figure}[b!]
\centering
\includegraphics[width=0.2\textwidth]{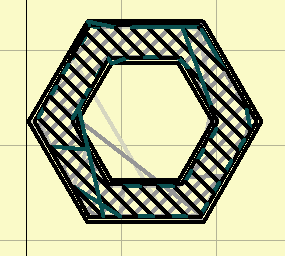}
\caption{An example of polygon with outer, inner, and filling parts.}
\label{fig:30gclo}
\end{figure}

{\bf Clustered 3DP.}
Let us first clarify that each layer is divided into different subpolygons and each subpolygon can be formed by at most three clusters of compulsory edges: the outer part (the border of the subpolygon), the inner part (that follows the outer part internally, to strengthen the structure), and the filling part (that represents the edges used to fill the subpolygon), see Figure \ref{fig:30gclo}.
The second algorithm we present is called {\em Clustered 3DP} algorithm, and follows a similar idea with respect to the Closest 3DP algorithm, while making use of the different clusters just presented.
It selects a cluster on a layer (choosing among outer, inner, and filling) and travels all the compulsory edges of the cluster, before going to another cluster. It starts from vertex 0 and selects the closest vertex of the compulsory edges on layer $L_1$. This action not only defines the next edge, but also the first cluster to be completed before moving to the next one. When the current cluster is completed, the algorithm moves to the closest vertex, which defines the next edge and the next cluster. The algorithm continues until all clusters of the current layer are completed and then the closest vertex of the next layer is chosen. When the last layer is completed we travel back to vertex 0 and the algorithm terminates.

{\bf Look Ahead 3DP.}
The name of the third algorithm that we propose is inspired by the fact that it makes use of heuristic information on the future possible decisions.
It starts from vertex 0 and then selects the three closest vertices among those of the compulsory edges of layer $L_1$.  A particular cost is then computed by running the Closest 3DP algorithm from each of the three vertices till the end of the layer. The cost is computed as follows: we sum up the traveling cost of the non-printing edges needed to terminate the layer by the Closest 3DP, the vertex of the three selected ones that furnishes with the cheapest cost is used. This process is iterated until all  compulsory edges of the current layer are used at least once. We then use the cheapest edge to move to the next layer and repeat the procedure until the last layer. To conclude the algorithm, we move back to the vertex 0.

{\bf Shortest Path Based 3DP.}
The forth proposed algorithm is called {\em Shortest Path Based 3DP} and is divided into two phases, the first one builds a supporting graph to be used in the second phase. The second phase computes the Dijkstra algorithm on the newly built supporting graph to obtain a solution for the 3DRP.
In the first phase we define the supporting graph $G'=(V',E')$, where $V' \subseteq V$ is the set of vertices and $E'$ the set of edges that we describe in the following.
For each layer we select a random set of vertices among those of the compulsory edges ($R$), called {\it starting vertices}, from which we can start the layer. The union between the starting vertices, vertex 0, and its duplication, vertex 0', defines $V'$. The set $E'$ is defined by: the edges going from each starting vertex of one layer to each starting vertex of the next layer, the edges going from vertex 0 to the starting point of the first layer $L_1$, and the edges going from the starting points of the last layer $L_{|\mathcal{L}|}$ to vertex 0'. Starting from each starting vertex we terminate the layer as if using the Closest 3DP algorithm ending up in one and only one vertex that we call {\it ending vertex}. We call $c'_{ij}$ the cost of each edge of $E'$. The costs $c'_{ij}$ is computed as the summation of the distance between a starting vertex and the corresponding ending vertex, whose cost is computed by the Closest 3DP algorithm, and the distance between the corresponding ending vertices and the starting vertices of the next layer. The starting vertices of the first layer are linked to the vertex 0 and a cost $c'_{0j} = c_{0j}, j \in L_1$  is associated to them. The ending vertices of the last layer are linked to the vertex 0', that has the same coordinates of vertex 0, and a cost $c'_{j0'}=c_{j0}, j \in L_{|\mathcal{L}|}$ is associated to them.  The costs $c'_{ij}, (i,j) \in V'$ give the weight of the edges of the supporting graph $G'$. In the second phase, we compute the shortest path form 0 to 0' on the supporting graph $G'$. This leads to a feasible solution for the 3DRP.

\section{Computational Results} \label{sec:comp}

We coded our algorithms in C++ and tested them on an Intel Core i3-2100 with 3.10 GHZ. We collected 30 instances on the Internet from the 3DP design community \url{thingiverse.com}. The instances have been processed by our algorithms and with one of the most famous software used for slicing polygonal meshes and for path definition, Cura \cite{cura}. We used the instances as input for Cura that takes the polygonal meshes of the instances, slice them, and defines a path of the nozzle in a code called g-code. The g-code represents a heuristic solution for the related 3DRP in a language that the 3D printer can understand. By processing the obtained g-code we can obtain the set of layers, the set of compulsory edges to be printed, the corresponding vertices and their coordinates. These pieces of information define the data that we used as instances for our algorithms. In Table \ref{tab:maintable}, we report the name of the instance, the number $|\mathcal{L}|$ of layers, the number $|R|$ of compulsory edges, and the total distance given by Cura.

Let  $\Gamma$ be the cost of the non-printing edges used in a 3DRP solution. Note that the compulsory edges are the same for all the solutions, and thus we report only the non-printing cost. Let $\Gamma_{Cura}$ be the value of $\Gamma$ proposed by Cura, in $mm$, and $\Gamma_{A}$ the value $\Gamma$ defined by the solution obtained thanks to one of our algorithms.  Let $\%impr$ be the percentage improvement between them, computed as $\%impr=100*(\Gamma_{A}-\Gamma_{Cura})/\Gamma_{Cura}$. In Table \ref{tab:maintable}, we report the results of the proposed Closest 3DP, Clustered 3DP, and Look Ahead 3DP algorithms showing their percentage improvement  and their computing time in seconds.

From Table \ref{tab:maintable}, one can observe that the Closest 3DP can improve the solution value for 28 out of 30 instances by 22.91\% on average. The average solving time is very promising being lower than one second. With Clustered 3DP algorithm we can improve the solution value with respect to Cura for 27 out of 30 instances, with an average improvement of 12.15\% in a time that is about half a second on average.
The Look Ahead (LA) 3DP algorithm results show that the solution value is improved for every instance with respect to Cura, by 36.73\% on average. To reach these results, the LA needs to run for about 39 minutes on average, which is a relevant solving time if compared with the previously proposed algorithms. These very good solutions pay their quality in terms of higher solving times.

In Table \ref{tab:shortest}, the average results for the Shortest Path Based (SPB) 3DP algorithms are shown, averaged on all the instances. In the left part, we report the results for the SPB, where the number of starting vertices (SV) for each layer is fixed to 10, 20, 30, 50, 60, and 100; in the right part, we present the results for the SPB where the SV are fixed with respect to a percentage of the vertices of each layer (1\%, 2\%, 5\%, 10\%, 15\%, 20\%, and 25\%).
One can see that for both SPB types the average improvement increases with the number of starting vertices explored. The improvement is proportional also to the computing times. The SPB with fixed numbers of starting vertices can reach a 27.38\% gap on average in the best case, that is represented by the starting vertices set to 100. This value is obtained in about 80 seconds. One can see that this version works similarly to the version with percentage starting vertices when evaluating the average improvement, but needs consistently more seconds to terminate.
The algorithm  that shows the best results in terms of improvement of the solution value is LA. However, LA needs higher computing times with respect to the other algorithms. The SPB furnishes, in its two versions, good results in acceptable solving times. %Similar but weaker improvements can be obtained by the Closest 3DP algorithm showing very small computing times.
\begin{table}[htbp]
  \centering
  \small
  \caption{Computational results for Closest, Clustered, and Look Ahead}
    \begin{tabular}{lrrrrrrrrr}
    \toprule
     \multicolumn{3}{c}{Instance} & \multicolumn{1}{c}{Cura}  & \multicolumn{2}{c}{Closest} & \multicolumn{2}{c}{Clustered} & \multicolumn{2}{c}{Look Ahead} \\
     \cmidrule(lr){1-3} \cmidrule(lr){4-4} \cmidrule(lr){5-6} \cmidrule(lr){7-8} \cmidrule(lr){9-10}
     Name  & \multicolumn{1}{c}{$|\mathcal{L}|$ } & \multicolumn{1}{c}{$|R|$} & \multicolumn{1}{c}{$\Gamma$ }  & \multicolumn{1}{c}{$\%impr$} & \multicolumn{1}{c}{$t(sec)$} & \multicolumn{1}{c}{$\%impr$} & \multicolumn{1}{c}{$t(sec)$} & \multicolumn{1}{c}{$\%impr$} & \multicolumn{1}{c}{$t(sec)$} \\
    \cmidrule(lr){1-3} \cmidrule(lr){4-4} \cmidrule(lr){5-6} \cmidrule(lr){7-8} \cmidrule(lr){9-10}
    30grams & 219 & 8693 & 11687.56 & 5.66 & 0.00 & 2.84 & 0.00 & 17.38 & 0.31 \\
    angel & 153 & 353671 & 147351.1 & 10.12 & 6.67 & 8.51 & 3.56 & 29.02 & 22664.21 \\
    block\_4x4 & 37 & 28694 & 15263.46 & 40.68 & 0.14 & 19.86 & 0.09 & 48.68 & 122.82 \\
    brickb & 332 & 107355 & 132462.2 & -1.00 & 0.59 & -1.03 & 0.47 & 12.67 & 765.48 \\
    bunny & 41 & 21333 & 7569.886 & 16.88 & 0.08 & 9.68 & 0.05 & 31.25 & 45.67 \\
    chamfer & 20 & 12700 & 11813.21 & 3.17 & 0.09 & -3.12 & 0.06 & 15.39 & 77.05 \\
    Cityscape & 342 & 488409 & 293773.5 & 9.95 & 7.44 & 5.02 & 3.64 & 30.59 & 17105.26 \\
    gear\_b & 85 & 19275 & 6167.369 & 3.86 & 0.03 & 7.19 & 0.03 & 14.26 & 10.42 \\
    geoboard & 82 & 52244 & 51043.13 & 3.93 & 0.28 & 2.89 & 0.24 & 15.14 & 303.95 \\
    hilbert2 & 99 & 26703 & 14938.52 & 11.71 & 0.03 & 11.71 & 0.08 & 21.50 & 13.39 \\
    Hilton\_Chicago & 73 & 7179 & 6081.438 & 18.09 & 0.00 & 14.28 & 0.00 & 30.21 & 1.06 \\
    ingranaggio & 39 & 52569 & 16516.57 & 22.22 & 0.42 & 4.34 & 0.33 & 36.70 & 758.60 \\
    miwin1 & 132 & 130267 & 95751.26 & 13.26 & 1.11 & 0.90 & 0.58 & 32.22 & 1645.93 \\
    orso\_bis & 52 & 35772 & 41535.23 & 60.24 & 0.16 & 38.17 & 0.11 & 73.46 & 131.84 \\
    pesce\_bis & 52 & 65262 & 77239.71 & 64.16 & 0.53 & 31.46 & 0.36 & 70.89 & 1012.54 \\
    polysoup & 3 & 3018 & 2042.693 & 71.14 & 0.02 & 43.72 & 0.02 & 78.06 & 21.67 \\
    portab & 32 & 47293 & 34371.95 & 13.53 & 0.44 & 12.01 & 0.45 & 32.08 & 1208.38 \\
    portaingr & 51 & 69664 & 17318.14 & 10.54 & 0.92 & 7.57 & 0.49 & 29.48 & 3032.91 \\
    RFin & 188 & 35898 & 9401.675 & 10.32 & 0.05 & 0.42 & 0.03 & 30.41 & 17.16 \\
    RostockBottom & 44 & 18120 & 27514.87 & 31.54 & 0.13 & 13.92 & 0.14 & 46.74 & 128.41 \\
    RostockTop & 50 & 21132 & 22071.05 & 23.36 & 0.17 & 10.80 & 0.13 & 36.72 & 369.54 \\
    Rpenta & 369 & 433954 & 203950.3 & 28.52 & 3.13 & 18.71 & 2.05 & 46.30 & 5144.17 \\
    RShowerHead & 274 & 146602 & 140072.8 & 24.85 & 0.52 & 17.41 & 0.45 & 40.69 & 472.31 \\
    tardisflat2 & 170 & 24695 & 23625.85 & 39.63 & 0.03 & 9.14 & 0.03 & 47.73 & 6.00 \\
    tetra & 160 & 20839 & 13115.34 & -4.81 & 0.03 & -0.89 & 0.03 & 16.62 & 6.94 \\
    trail & 125 & 268372 & 244951.1 & 48.29 & 4.19 & 45.74 & 2.77 & 59.47 & 14259.08 \\
    wbuilding & 100 & 3364 & 2643.942 & 2.86 & 0.00 & 1.21 & 0.02 & 16.70 & 0.08 \\
    wine\_fixed & 39 & 20276 & 29576.97 & 38.03 & 0.09 & 12.84 & 0.06 & 48.10 & 62.18 \\
    wine1 & 39 & 19145 & 29808.14 & 40.19 & 0.08 & 1.59 & 0.06 & 51.12 & 46.33 \\
    WitchCastle & 112 & 56334 & 12736.63 & 26.33 & 0.22 & 17.49 & 0.17 & 42.37 & 187.03 \\
    \cmidrule(lr){1-3} \cmidrule(lr){4-4} \cmidrule(lr){5-6} \cmidrule(lr){7-8} \cmidrule(lr){9-10}
    Average &   &   &   & 22.91 & 0.92 & 12.15 & 0.55 & 36.73 & 2320.69 \\
        \bottomrule
    \end{tabular}%
  \label{tab:maintable}%
\end{table}%
\begin{table}[htbp]
  \centering
  \caption{Results for the Shortest Path Based 3DP algorithm.}
    \begin{tabular}{rrrlrr}
          \toprule
    \multicolumn{1}{l}{$\#SV$} & \multicolumn{1}{l}{$\%impr$} & \multicolumn{1}{l}{$t(sec)$} & $\%SV$ & \multicolumn{1}{l}{$\%impr$} & \multicolumn{1}{l}{$t(sec)$} \\
         \cmidrule(lr){1-3} \cmidrule(lr){4-6}
    10 & 23.50 & 8.59 & 1\% & 18.54 & 24.96 \\
    20 & 25.20 & 16.32 & 2\% & 21.10 & 48.74 \\
    30 & 26.04 & 24.16 & 5\% & 24.02 & 118.13 \\
    40 & 26.43 & 31.85 & 10\% & 25.70 & 181.53 \\
    50 & 26.63 & 39.83 & 15\% & 26.29 & 200.25 \\
    60 & 26.89 & 47.51 & 20\% & 26.59 & 208.47 \\
    100 & 27.38 & 80.94 & 25\% & 27.00 & 212.81 \\
        \bottomrule
    \end{tabular}%
  \label{tab:shortest}%
\end{table}%

\section{Conclusions and Future Research Directions} \label{future}
We studied a generalization of the RPP related to the definition of the tool path in a 3DP process, called 3DRP.
We provided a collection of 3DRP instances and solved them with four heuristic algorithms. We compared them with one of the most famous software for 3DP path definition, Cura \cite{cura}. We obtained very good results, consistently improving the Cura solution values within limited computing times. The Closest and the Clustered algorithms can improve the solutions in less than a second for most of the instances; the Look Ahead algorithm can improve the solutions by a 36\% but requires long times; on the other hand we believe that the Shortest Path Based algorithms provides the best trade off between improvements and computing times.
The obtained data are a simulation of the real print, in the future we intend to test the algorithms by performing real prints of all the objects.

Other versions of the 3DRP with other secondary features can be considered. In the following we report an overview of these features.
%Some of these will be part of next studies.\\

\noindent {\bf Quality constraints.} The 3DP process can produce some blemish on the external surface of the printed  product. Indeed, when traveling the optional edges after ejecting material, the nozzle can eject some leftovers and ruin the quality of the external surface. Additional constraints could be considered to avoid it.

%\begin{enumerate}
%\item
{\it Second level of clusterization: Subpolygons.} Several separated polygons, called {\em subpolygons}, can be present in the same layer. Subpolygons produce a second level of clusterization other than the layers. Defining the nozzle path by following this second level of clusterization could improve the quality of the external surface. To diminish blemish on the external surface, we can impose a new constraint: to complete a polygon before stepping to the next one, by taking one and only one optional arc between a subpolygon to another on the same layer.

{\it Corner.} To preserve the quality of the surface from the blemishes, we can consider not only to leave the polygon just once, but also to impose the nozzle to pass only through a vertex of the outer part of the same subpolygon, called {\em corner}.
%, in order to avoid the intersection of the tool path and the surface when this has already been built.
By doing this, we will leave the majority of the blemishes on a corner and make it easier to remove them after the print. %To account for this we would have to solve a problem called {\it Shortest Path with Linear Obstacles.}

%\item
{\it Temperature.} The 3DP process includes the melting of the material before extruding it from the nozzle. After being melt, the material conserves a high temperature for a certain amount of time. Thus material cannot be extruded on the top of the previous layer until its temperature is not under a certain value, otherwise flaws or bending can happen. We can impose this by setting that an edge cannot be completed before an underlying edge had the time to cool down.
%Since an edge can be not completely overlapping with respect to those of the next/previous layer we must define a set of edges $C(i)$ of the layer $L_h$, that must be completed not before a given amount of time after visiting the edge $i$ of layer $L_{h-1}$ has passed.
%In order to cool down the temperature of some edges, the speed of the fan can also be considered as a variable.

%\item
{\it First vertex of a layer.} When moving from a non-printing edge to a printing one a small excess of material can be extruded. This can lead to an accumulation of unwilled material on the surface, especially when the first vertex of a layer has the similar x-y coordinates with respect to the first vertices of the previous layers. To avoid this, we can impose that the first vertex of a layer must have different coordinates if compared to the first vertex of the previous layer.
%Indeed, the same line of reasoning can be used when starting overlapping subpolygons on different layers, and not just the first vertex.

%{\it Quality constraints.} We can consider some features in order to avoid some superficial flaws:
%	 \begin{itemize}
%		\item Moving inside the subpolygon so to avoid to cross the crust if it is not necessary.
%		\item When moving from a subpolygon to the next one: only passing throw a corner of the subpolygon.
%		\item When starting the next layer: not to start from the same starting point of the previous one.
%	\end{itemize}
%\begin{enumerate}
%\item {\it Typical features of a RPP and general characteristics:}
%	\begin{itemize}
%		\item The solution must start from the depot, and terminate in the same vertex (the depot).
%		\item The solution must visit at least once the compulsory edges/arcs.
%		\item It is also allowed to use the optional edges/arcs -as many times as needed- to obtain a eulerian cycle.
%		\item Connectivity of the solution must be imposed.
%	\end{itemize}
%\end{enumerate}

{\bf Costs.}
We considered the cost of the edges as the geometric distance between two vertices and the time to travel the edge as directly linear to the distance. We could, alternatively, model them with an acceleration at first, a constant speed in the middle, and a deceleration at the end. Moreover, the deceleration could depend from the angle between the current edge and the following one visited, and the acceleration from the angle between the current edge and the previous one. This can lead to different modeling options, considering non-linear distances between the vertices or considering the speed of the nozzle as a variable.
%Taking into account the acceleration and deceleration needed to obtain the maximum possible speed on a given edge, we can recompute the cost of each edge starting from the distance. If we considering the angles between the previous and following edges, then the problem must be considered as one with directed arcs instead of edges. .
%\item{\it Clustered problem: Layers.}
%	\begin{itemize}
%	\item Between a layer and the subsequent one exist only optional edges/arcs and one and only one can be used.
%	\item Between two non subsequent layers no edge/arc is present/used.
%	\item A layer cannot be started until the previous one has not been completed.
%	\end{itemize}

{\bf Other Features.}
%\begin{enumerate}
%\item
{\it Layers with different thickness.} Sometimes, to deal with structural problems or just to avoid to obtain a texture too subtle inside the printed object, we can print some parts of the same layer with a different thickness. %Two different thicknesses must be related (one must be a multiple of the other). %This can be handled considering diverse relationships between layers and clusters.
%\item
{\it Supporting material.} Some 3D objects can have undercuts, thus a supporting material is needed when printing. In order to minimize the supporting material, the positioning of the object can be a variable to account for.
%\item
{\it Filling.} In this work we consider the filling part of a certain polygon as an input, but to decide the filling, in terms of distances between the edges of the filling, their direction, and even a particular pattern can change the optimization of the print or be part of it.
%{\it Mesh.} Triangulation/Meshes: how to do it? For us it is an input.
%\item
{\it Colors.} Some printing objects can have different set of materials and/or colors to be used. Indeed we can see the supporting material part of this case. When printing an object made of different colors or materials we can decide to print first the polygons of one material and then the others by deciding a sequence. Otherwise, we can decide to mix the use of colors to make the printing time shorter. Then the number of nozzles becomes relevant: if we have more materials than nozzles we must consider the time spent in changing material and going back to the vertex 0 to clean the nozzle, otherwise we must consider the distance between the nozzle of one color and the others.
%\item{\it Buffering.} We can try to buffer the algorithm printing the solution and the g-code feed to the 3D printer.
%\end{enumerate}

All these features have a strong impact on the speed of the print and on the quality of the surface, both very important characteristics in 3DP. We believe  that they will be the focus of many future research works in optimization.

\end{document}